\documentclass[11pt]{amsart}

\usepackage{amsmath}
\usepackage{amssymb}
\usepackage{amsthm}
\usepackage{xcolor}
\usepackage{array}
\usepackage[margin=1.15in]{geometry}
\usepackage[colorlinks=true,linkcolor=blue,urlcolor=blue,citecolor=blue]{hyperref}
\usepackage[capitalise]{cleveref}

\newcommand{\col}[2]{\left(\begin{array}{@{\,}#1@{\,}}#2\end{array}\right)}
\theoremstyle{plain}
\newtheorem{theorem}{Theorem}[section]
\newtheorem{lemma}[theorem]{Lemma}
\newtheorem{corollary}[theorem]{Corollary}

\theoremstyle{definition}
\newtheorem{definition}[theorem]{Definition}
\newtheorem{example}[theorem]{Example}

\theoremstyle{remark}

\DeclareMathOperator{\std}{std}
\DeclareMathOperator{\inv}{inv}
\DeclareMathOperator{\maj}{maj}

\DeclareMathOperator{\foata}{foata}
\DeclareMathOperator{\spot}{spot}
\DeclareMathOperator{\plate}{plate}
\DeclareMathOperator{\blk}{block}
\DeclareMathOperator{\pref}{pref}
\newcommand{\UPF}{\mathrm{UPF}}
\newcommand{\A}{\mathcal{A}}
\newcommand{\rd}[1]{\textcolor{red}{#1}}
\newcommand{\bl}[1]{\textcolor{blue}{#1}}

\title{Foata and Spot Commutation for Parking Functions}
\author{Eshaan Sombhatta}
\address{Thomas Jefferson High School for Science and Technology}
\email{2027esombhat@tjhsst.edu}
\subjclass[2020]{05A05, 05A15, 05A19}
\keywords{Parking function, unit interval parking function, Foata bijection, inversion, major index}

\begin{document}

\begin{abstract}

For a parking function $\alpha$, let $\spot(\alpha)$ be the permutation whose $i$th entry is the space occupied by car $i$ and let $\foata(w)$ be Foata's second fundamental transformation on words. 
Conjecture 6.3 of Celano et al. asserts that $\foata(\spot(\alpha)) = \spot(\foata(\alpha))$ for every unit interval parking function $\alpha$.
This paper proves the identity for a larger class: parking functions where each prime component (borrowing Gessel's definition) of $\alpha$ is nondecreasing. 
Since every unit interval parking function is composed of nondecreasing prime parking functions, the conjecture follows. 
A key step is showing that $\foata$ preserves the set of parking functions $\alpha$ with $\spot(\alpha) = \std(\alpha)$, where $\std$ denotes standardization.
In particular, this proves that inversion count and major index are equidistributed over this larger class.
\end{abstract}

\maketitle

\section{Introduction}

Parking functions, introduced by Konheim and Weiss \cite{KW}, come from the following problem. There are $n$ cars parking on a one-way street with $n$ parking spaces labeled $1$ through $n$. The $i$th car prefers spot $\alpha_i \in \{1,  \dots, n \}$ and parks there if it is unoccupied. Otherwise, it continues to the right until it finds an available parking space. If any car is unable to park, the parking process is unsuccessful. If all cars park successfully, then the tuple \[ \alpha=(\alpha_1,\alpha_2,\ldots,\alpha_n) \] is called a \textbf{parking function}. $\mathrm{PF}_n$ denotes the set of parking functions of length $n$. 
The \textbf{spot permutation} $\spot(\alpha)$ of a parking function $\alpha$ is the permutation whose $i$th entry is the spot in which the $i$th car parked. 

Though rearranging the cars of a parking function produces another parking function, it can change the spot permutation. This motivates the study of how permutations of a parking function affect the parking process. 

A natural rearrangement to consider is \textit{Foata's second fundamental transformation} on words, a content-preserving
bijection satisfying
\[
    \operatorname{inv}(\operatorname{foata}(w))
    = \operatorname{maj}(w),
\]
where $\operatorname{inv}$ is the \textit{inversion count} and $\operatorname{maj}$ is the \textit{major index} \cite{Foata}. Since $\foata$ preserves content, it sends $\mathrm{PF}_n$ to $\mathrm{PF}_n$.

This paper builds on the work of Celano, Elder, Hadaway, Harris, Martin,
Priestley, and Udell \cite{CEHHMPU}. Their Conjecture~6.3 asserts that for every \textbf{unit interval parking function} $\alpha$, in which each car parks
at most one space beyond its preference,
\[ \operatorname{foata}(\operatorname{spot}(\alpha)) = \operatorname{spot}(\operatorname{foata}(\alpha)). \]

Define the \textbf{standardization} of a word to be the permutation $\std(w)$ that ranks the entries of $w$ from smallest to largest, breaking ties by position, and define the set of \textbf{standard parking functions} as \[ \A_n = \{ \alpha\in\mathrm{PF}_n : \operatorname{std} (\alpha)=\operatorname{spot}(\alpha) \}. \]
The main result of this paper is that $\foata(\spot(\alpha)) = \spot(\foata(\alpha))$ holds whenever $\alpha\in\A_n$. This paper characterizes $\A_n$ as the set of parking functions whose \textit{prime blocks} are weakly increasing. This generalizes the block characterization of unit interval parking functions by Bradt, Elder, Harris, Rojas Kirby, Reutercrona, Wang, and Whidden \cite{Bradt} by replacing the condition that each prime block be a prime unit interval parking function with the condition that each prime block be a prime weakly increasing parking function.
Since the entries in blocks of unit interval parking functions are weakly increasing, $\A_n$ contains all unit interval parking functions. For $n\geq3$, this containment is strict, so the result establishes Conjecture 6.3 for a larger class.

Membership in $\A_n$ is sufficient but not necessary for the
identity. For example, $1211\notin\A_4$, but
\[
    \operatorname{foata}(\operatorname{spot}(1211))
    = \operatorname{spot}(\operatorname{foata}(1211))
    = 1234.
\]
A full characterization of the parking functions satisfying the identity remains open.

\section{Preliminaries}
\begin{definition}
The standardization $\std(w)$ of a word $w = (w_1,\dots,w_n)$ is the
permutation of $1$ through $n$ that ranks the entries of $w$ from smallest to
largest, breaking ties by position (earlier entries rank lower).
\end{definition}

\begin{definition}
Let $\A$ be the set of parking functions $\alpha$ with
$\std(\alpha) = \spot(\alpha)$.
Write $\A_n = \A \cap \mathrm{PF}_n$.
\end{definition}

\begin{definition}
    Let $\UPF$ be the set of parking functions $\alpha$ where $\max_i(\spot(\alpha)_i - \alpha_i) \leq 1$. Therefore, $\UPF_n$ is the set of parking functions of length $n$ where the maximum displacement is at most $1$.
\end{definition}

\begin{definition}
An integer $b$ with $1\leq b\leq n$ is a \textbf{breakpoint} of a parking function $\alpha\in\mathrm{PF}_n$ if exactly $b$ cars prefer the first $b$ spots. A parking function is \textbf{prime} if its only breakpoint is $n$. For a discussion of prime parking functions, see \cite[Section~2]{Meyles}.
    
Throughout this paper, the \textbf{prime blocks} of $\alpha$ are the subsequences with preferences in $[a+1,b]$, listed in the same order in which the corresponding cars arrive in $\alpha$, where $a<b$ are consecutive breakpoints. Subtracting $a$ from every preference in such a block shifts its preference interval from $[a+1,b]$ to $[1,b-a]$; the resulting word is a prime parking function \cite{GesselSeo}. For a car $i$, write $\blk(i)$ for the unique prime block containing its preference $\alpha_i$.
\end{definition}

\begin{definition}
Following \cite{Foata}, \textbf{Foata's second fundamental transformation} is the function defined inductively on words $w=w_1\cdots w_n$ as follows.
For $n=1$, set $\foata(w)=w$.

For $n>1$, $\foata(w)$ is defined as follows.  Let $w'=\foata(w_1\cdots w_{n-1})=w'_1\cdots w'_{n-1}$. 
If $w_n\geq w'_{n-1}$, then place a separator after every $w'_i$ such that $w_n\geq w'_i$.
If $w_n<w'_{n-1}$, then place a separator after every $w'_i$ such that $w_n<w'_i$.
Cycle each \textbf{segment} (each maximal subword between consecutive separators) by moving its rightmost entry to the start.
Finally, append $w_n$.  The result is $\foata(w)$. 

A step-by-step illustration of this procedure is given in \cref{ex:4.1}.
\end{definition}

\begin{definition}
    For a word $w$, an inversion is a pair of indices $(i, j)$, where $i < j$, but $w_i > w_j$. The \textbf{inversion count} of a word $w$, denoted $\inv(w)$, is the number of inversions in a word $w$.
\end{definition}
\begin{definition}
    The \textbf{major index} of a word, denoted $\maj(w)$, is the sum of the indices $i$ where $w_i > w_{i+1}$.
\end{definition}
\section{Proof structure}

For $\alpha\in\A$, the proof establishes the following chain of equalities:
\[
  \foata(\spot(\alpha))
    \underset{\spot \,=\, \std}{=} \foata(\std(\alpha))
    \underset{\foata \,\circ\, \std \,=\, \std \,\circ\, \foata}{=}
      \std(\foata(\alpha))
    \underset{\foata(\alpha) \,\in\, \A}{=} \spot(\foata(\alpha))
\]

The first equality, $\foata(\spot(\alpha)) = \foata(\std(\alpha))$ holds since for all $\alpha$ in $\A$, $\spot(\alpha) = \std(\alpha)$. 
The second equality, $\foata(\std(\alpha)) = \std(\foata(\alpha))$ is an identity which holds more generally for all words $w$, and will be proven in \cref{sec:4}.
The third equality, $\std(\foata(\alpha)) = \spot(\foata(\alpha))$ is a consequence of the fact that $\foata$ preserves membership in $\A$. The statement $\foata(\alpha) \in \A$ will be proven in \cref{sec:5}.

\section{Foata commutes with standardization}
\label{sec:4}

\begin{example}
\label{ex:4.1}
The dividers occur at the same positions for $\foata(\std(w))$ and $\foata(w)$.  Below is an illustrative example.
\begin{center}
\begin{tabular}{|l|l|}
\hline
$w = 314132$ & $\std(w) = 416253$ \\
\hline
$3 \to 3$                   & $4 \to 4$ \\
$3\mid 1 \to 31$            & $4\mid 1 \to 41$ \\
$3\mid 1\mid 4 \to 314$     & $4\mid 1\mid 6 \to 416$ \\
$3\mid 14\mid 1 \to 3411$   & $4\mid 16\mid 2 \to 4612$ \\
$3\mid 41\mid 1\mid 3 \to 31413$ & $4\mid 61\mid 2\mid 5 \to 41625$ \\
$3\mid 14\mid 13\mid 2 \to 341312$ & $4\mid 16\mid 25\mid 3 \to 461523$ \\
\hline
\end{tabular}
\end{center}
\end{example}

To formalize this observation, label the letters by their original positions.

\begin{definition}[Two-line notation]
Let $w = w_1\cdots w_n$ be a word, and write it as the two-line array
\[
  \begin{pmatrix}
    1 & 2 & \cdots & n \\
    w_1 & w_2 & \cdots & w_n
  \end{pmatrix},
\]
so each letter is labeled by its position. 

For words $a$ and $b$, where $b$ is a rearrangement of $a$, $\plate_b (a)$ is the permutation obtained by labeling the entries of $b$ by their positions, rearranging these labeled entries to form $a$ while preserving the relative order of equal entries, and reading off the labels.

When applying $\foata$ to a two-line word, apply $\foata$ to the columns, comparing only the bottom entries, but cycling both the top and bottom entry. Then $\plate_w(\foata(w))$ is the top row of the result. 
\end{definition}

\begin{example}
    For example, for $w = 314132$, we have $\plate_w(\foata(w)) = 132546$. Below, $\foata$ is applied step-by-step to the two-line array.

\begin{center}
\renewcommand{\arraystretch}{1.15}
\begin{tabular}{|c|l|}
\hline
Foata step & $\foata$ applied to the two-line array \\
\hline
$1$ & $\col{c}{1\\3} \;\to\; \col{c}{1\\3}$ \rule[-3ex]{0pt}{7ex}\\
$2$ & $\col{c|c}{1&2\\3&1} \;\to\; \col{cc}{1&2\\3&1}$ \rule[-3ex]{0pt}{7ex}\\
$3$ & $\col{c|c|c}{1&2&3\\3&1&4} \;\to\; \col{ccc}{1&2&3\\3&1&4}$ \rule[-3ex]{0pt}{7ex}\\
$4$ & $\col{c|cc|c}{1&2&3&4\\3&1&4&1} \;\to\; \col{cccc}{1&3&2&4\\3&4&1&1}$ \rule[-3ex]{0pt}{7ex}\\
$5$ & $\col{c|cc|c|c}{1&3&2&4&5\\3&4&1&1&3} \;\to\; \col{ccccc}{1&2&3&4&5\\3&1&4&1&3}$ \rule[-3ex]{0pt}{7ex}\\
$6$ & $\col{c|cc|cc|c}{1&2&3&4&5&6\\3&1&4&1&3&2} \;\to\; \col{cccccc}{1&3&2&5&4&6\\3&4&1&3&1&2}$ \rule[-3ex]{0pt}{7ex}\\
\hline
\end{tabular}
\end{center}

\end{example}

The following lemma shows the equality of the resulting label orders.

\begin{lemma}
For all words $w$, $\plate_w(\foata(w)) = \plate_{\std(w)}(\foata(\std(w)))$.
\end{lemma}

\begin{proof}
Proceed by induction on prefix length. The base case of one letter is clear. Now, assume that both runs have the same label order before appending $w_q$. 
Since $\std$ acts positionwise, car $i$ of $w$ is car $i$ of $\std(w)$. The last letters therefore choose the same divider rule since if $p < q$, then $w_p \leq w_q \iff \std(w)_p < \std(w)_q$.

At the step appending $w_q$, every divider decision is a test of some
        earlier letter $w_p$, $p < q$, against $w_q$.
If $w_p \neq w_q$, then
        $w_p < w_q \iff \std(w)_p < \std(w)_q$.  So the divider placement is
        identical.

If $w_p = w_q$, then $\std(w)_p < \std(w)_q$, since $p < q$.
First suppose the condition is that a divider is placed after all
            elements $e$ satisfying $e \le w_q$.  Then in $w$, a divider will be
            placed after $w_p$ since $w_p = w_q$, and in $\std(w)$, a divider will be placed
            since $\std(w)_p < \std(w)_q$.
Now suppose the condition is that a divider is placed after all
            elements $e$ satisfying $e > w_q$.  In $w$, no divider will be
            placed, since $w_p = w_q$, which is not greater than $w_q$.  In
            $\std(w)$, no divider will be placed, since
            $\std(w)_p < \std(w)_q$, which is not greater than $\std(w)_q$.

So both runs cut into the same compartments at every step.  Therefore,
        both runs apply the same permutation.
\end{proof}

\begin{corollary}
$\foata(\std(w)) = \std(\foata(w))$.
\end{corollary}

\begin{proof}
Write $\sigma = \plate_w(\foata(w))$, so that $\foata(w) = w \circ \sigma$ and
$\foata(\std(w)) = \std(w) \circ \sigma$.  It suffices to prove that $\sigma$ preserves the
relative order of the occurrences of each letter, implying
\[
  \foata(\std(w)) = \std(w) \circ \sigma
    \underset{\sigma \text{ preserves order}}{=} \std(w \circ \sigma)
    = \std(\foata(w)).
\]
There is no crossing of equal letters during the $\foata$ process.
The ordering of two letters $w_p$ and $w_q$ can only be swapped if they are in the same compartment during the $\foata$ process.  Moreover, at the first swap of equal letters, if $p < q$, then $w_q$ must be the last letter in the compartment.

Suppose each element is being compared to an incoming letter $e$ in the
        $\foata$ process.
First, suppose that the condition is a divider is placed after all
            items greater than $e$.  Then $w_p = w_q$ but $w_q > e$ while $w_p$
            is not, a contradiction.
Next, suppose that the condition is a divider is placed after all
            items less than or equal to $e$.  Then $w_p = w_q$ but
            $w_q \le e$ while $w_p$ is not, a contradiction.
Finally, appending an item clearly cannot reverse equal occurrences. Therefore, $\std(w) \circ \sigma = \std(w \circ \sigma)$ and
$\foata(\std(w)) = \std(\foata(w))$.
\end{proof}

\section{Foata preserves membership in $\A$}
\label{sec:5}

The identity
$\std \circ \foata = \foata \circ \std$ holds for all words. The remaining argument concerns the parking process and uses the following notation.

\begin{definition}
Let $\spot(\alpha_i)$ denote the spot of the car initially at index $i$.
        For example, if $\alpha = 111$, then $\spot(\alpha_2) = 2$.
Let $\pref(\alpha)_i$ be the preference of the car that parks at spot
        $i$.
\end{definition}

The following examples illustrate membership in $\A$.

\begin{example}
The parking functions $112$, $131$, $123$, and $1122773388$ are in $\A$, while
$121$ is not.
\end{example}

The next lemma characterizes $\A$.

\begin{lemma}
$\alpha \in \A$ if and only if all prime blocks of $\alpha$ are in nondecreasing order. 
\end{lemma}

\begin{proof}
For sufficiency, note that $\spot(\alpha) = \std(\alpha)$ if $\alpha$ is nondecreasing, and prime blocks park independently.  
A permutation $\pi$ equals $\std(w)$ if
$\pi(i) < \pi(j)$ when $w_i < w_j$, and
$\pi(i) < \pi(j)$ when $w_i = w_j$ and $i < j$.
If all prime blocks of $\alpha$ are nondecreasing, this holds.
If $\blk(i) \neq \blk(j)$, then $\alpha_i < \alpha_j \iff \spot(\alpha_i) < \spot(\alpha_j)$ since blocks park independently.
If $\blk(i) = \blk(j)$, then whenever $\alpha_i < \alpha_j$, nondecreasing order gives $i < j$, and therefore $\spot(\alpha_i) < \spot(\alpha_j)$. If instead $\alpha_i = \alpha_j$, then $i < j \iff \spot(\alpha_i) < \spot(\alpha_j)$.

For necessity, assume $\alpha\in\A$.  It suffices to show that if $i < j$ but $\alpha_i > \alpha_j$, there is a breakpoint $k$ with $\alpha_j \leq k < \alpha_i$.

Choose $k$ to be the \emph{highest} spot below $\alpha_i$ still empty immediately after car $i$ parks. This set is nonempty, since the eventual spot of car $j$ is a member of it.
Since $\alpha \in \A$ and $\alpha_j < \alpha_i$,
        $\std(\alpha)_j < \std(\alpha)_i$ and
        $\spot(\alpha_j) < \spot(\alpha_i)$.
All the spots in $[\alpha_i, \spot(\alpha_i)]$ are occupied after car
        $i$ parks, so car $j$ must park before $\alpha_i$.

To show that $k$ is a breakpoint, it suffices to prove that there is no car with index $c$ where $\alpha_c \leq k$ but $\spot(\alpha_c) > k$.
If $c < i$, then space $k$ was empty when $c$ parked. Therefore $\spot(\alpha_c) < k$, since no car could park past the empty space at $k$. 
If $c > i$, then since $\alpha_c \leq k < \alpha_i$ and $\alpha \in \A$, it follows that $\spot(\alpha_c) < \spot(\alpha_i)$. After car $i$ parks, every space in $[k+1, \alpha_i - 1]$ is occupied by maximality of $k$. Further, every space in $[\alpha_i, \spot(\alpha_i)]$ is occupied. So car $c$ cannot park anywhere in $[k + 1, \spot(\alpha_i)]$ and $\spot(\alpha_c) \leq k$.
$c \neq i$ since $\alpha_i > k$.
Therefore, $k$ is a breakpoint and $i$ and $j$ are in different blocks.
\end{proof}

\begin{lemma}
For $\alpha\in\A$, $\foata$ preserves nondecreasing prime blocks.
\end{lemma}
Write $\foata^{(m)}(\alpha)=\foata(\alpha_1\cdots\alpha_m)$.
The superscript denotes the transformed prefix length, not iteration of the map.

\begin{example}
Consider $\alpha = 51516237$.  The Foata process proceeds as follows
(prime blocks are colored).
\begin{center} 
\begin{tabular}{|c|l|}
\hline
Foata step & $\foata^{(n)}(\alpha)$, \quad
             $\alpha = \bl{5}\rd{1}\bl{5}\rd{1}\bl{6}\rd{2}\rd{3}\bl{7}$ \\
\hline
$1$ & $\bl{5} \to \bl{5}$ \\
$2$ & $\bl{5}\mid\rd{1} \to \bl{5}\rd{1}$ \\
$3$ & $\bl{5}\mid\rd{1}\mid\bl{5} \to \bl{5}\rd{1}\bl{5}$ \\
$4$ & $\bl{5}\mid\rd{1}\bl{5}\mid\rd{1} \to \bl{55}\rd{11}$ \\
$5$ & $\bl{5}\mid\bl{5}\mid\rd{1}\mid\rd{1}\mid\bl{6}
       \to \bl{55}\rd{11}\bl{6}$ \\
$6$ & $\bl{5}\mid\bl{5}\mid\rd{11}\bl{6}\mid\rd{2}
       \to \bl{556}\rd{112}$ \\
$7$ & $\bl{556}\rd{1}\mid\rd{1}\mid\rd{2}\mid\rd{3}
       \to \rd{1}\bl{556}\rd{123}$ \\
$8$ & $\rd{1}\mid\bl{5}\mid\bl{5}\mid\bl{6}\mid\rd{1}\mid\rd{2}\mid\rd{3}\mid\bl{7}
       \to \rd{1}\bl{556}\rd{123}\bl{7}$ \\
\hline
\end{tabular}
\end{center}
Each prime block remains nondecreasing after each step.
\end{example}

\begin{proof}
Proceed by induction on the prefix length $n$, tracking the blocks of the original parking function $\alpha$. The claim holds for a one-letter prefix. Assume that the letters of each block occur in nondecreasing order in $\foata^{(n-1)}(\alpha)$.  Since $\alpha \in \A$, the letters
of each prime block arrive in nondecreasing order, so the incoming letter $\ell$
is at least every letter of its own prime block already placed, and appending to
the end maintains nondecreasing order.  It remains to check the cycling, which changes
the relative order between the last letter $x$ of a divider block and the
letters $y$ it is cycled backwards past.

Suppose the last letter of $\foata^{(n-1)}(\alpha)$ is $\le \ell$, so a
        divider is placed after each $e \le \ell$.  Every divider block then
        ends in a letter $\le \ell$ and contains no other such letter, so
        $x \le \ell < y$.  If $x$ and $y$ were in the same prime block, the
        invariant would give $y \le x$, since $y$ precedes $x$.  This is a contradiction.

Otherwise a divider is placed after each $e > \ell$, so
        $x > \ell \ge y$.  Suppose $x$ and $y$ lie in the same prime block $P$.
        Prime blocks correspond to intervals of preferences and $y \le \ell < x$, so
        $\ell \in P$.  But $x \in P$ is already placed and $x > \ell$,
        contradicting that $\ell$ is the largest letter of $P$ placed so far.
Either way $x$ crosses only letters of other prime blocks.
\end{proof}

\section{Main result and enumeration}

\begin{theorem}
For every $\alpha\in\A_n$,
\[
\foata(\spot(\alpha)) = \spot(\foata(\alpha)).
\]
In particular, the identity holds for every unit interval parking function.
\end{theorem}

\begin{proof}
The preceding results give
\[
\foata(\spot(\alpha)) = \foata(\std(\alpha))
= \std(\foata(\alpha)) = \spot(\foata(\alpha)).
\]
Further, $\A_n$ contains $\UPF_n$, since all prime blocks of every unit interval parking function are nondecreasing. Thus the theorem establishes Conjecture~6.3 of \cite{CEHHMPU}.
\end{proof}

Since standard parking functions can be decomposed into nondecreasing prime parking functions, and there are $C_{n-1}$ nondecreasing prime parking functions of length $n$, where $C_n$ is the $n$th Catalan number, we have the relation $|\A_0|= 1$ and 
\[ |\A_n| = \sum_{m = 1}^{n} {n \choose m}  C_{m-1}  |\A_{n-m}|. \]
Setting $A(z) = \sum_{n \geq 0} |\A_n| \frac{z^n}{n!}$ and $C(z) = \sum_{n \geq 0} C_n \frac{z^n}{n!}$, we have
\[ A(z) = 1 +  \left( \sum_{n \geq 1} C_{n-1} \frac{z^n}{n!} \right) \left( \sum_{n \geq 0} |\A_n| \frac{z^n}{n!} \right) = 1+  \left( \int_{0}^{z} C(t) \mathrm{d}t \right) A(z) \]
and 
\[A(z) = \frac{1}{1- \int_{0}^{z} C(t) \mathrm{d}t }.\]
The exponential generating function for the Catalan numbers is
\[ C(z)=e^{2z}\bigl(I_0(2z)-I_1(2z)\bigr), \]
where $I_n(z)$ is the $n$th ``modified Bessel function of the first kind.'' Using the identities
\[ \int_0^z e^tI_0(t)\,\mathrm dt = ze^z\bigl(I_0(z)-I_1(z)\bigr) \]
and
\[ \int_0^z e^tI_1(t)\,\mathrm dt = e^z\bigl((1-z)I_0(z)+zI_1(z)\bigr)-1\]
from \cite{marcum1948statistical}, we can write
\begin{align*}
A(z)
&= \frac{1}{1-\int_0^z C(t)\,\mathrm dt}\\
&= \frac{1}{1-\int_0^z e^{2t}\bigl(I_0(2t)-I_1(2t)\bigr)\,\mathrm dt}\\
&= \frac{1}{1-\frac12\int_0^{2z}e^t\bigl(I_0(t)-I_1(t)\bigr)\,\mathrm dt}\\
&= \frac{1}{1-\left(\frac12+\frac12e^{2z}\left[(4z-1)I_0(2z)-4zI_1(2z)\right]\right)}\\
&= \frac{2}{1+e^{2z}\left[(1-4z)I_0(2z)+4zI_1(2z)\right]}.
\end{align*}

Further, it is possible to find a $q$-exponential generating function $A(q,z)$. Using a nearly identical argument to the unit interval case of \cite{CelanoEtAl2025Inversions}, we find 
\[ A(q, z) = \frac{1}{1 - \sum_{n \geq 1} C_{n-1} \frac{z^n}{[n]_q!} }. \]

\section*{Acknowledgments}

The author thanks Gabe Udell for his guidance and comments on many drafts. The author also thanks Pamela Harris and Jennifer Elder for their helpful correspondence.

\bibliographystyle{amsalpha}
\bibliography{refs}

\end{document}